\documentclass[11pt]{article}
\usepackage{amsfonts,amsmath}
\usepackage{latexsym}
\usepackage{amsthm}
\usepackage{amscd}
\usepackage{commath}
\usepackage{epsfig}
\usepackage{amssymb}
\usepackage{tikz}
\usepackage{graphicx}
\usepackage{caption}
\usepackage{enumitem}
\usepackage{mathtools}
\usepackage[toc,page]{appendix}
\usepackage{hyperref}

\newcommand{\treeleaf}{\begin{tikzpicture}
		\node [fill=black, inner sep=2.5 pt,label=182: $v_0$] at (0,0) {} ;
		\draw[line width=0.05cm, black] (0,0) -- (0,1)  {};
		\node [fill=black, inner sep=2.5 pt,label=182: $v_1$] at (-1.8,3) {} ;
		\node [fill=black, inner sep=2.5 pt,label=182: $v_2$] at (-0.6,3) {} ;
		\node [fill=black, inner sep=2.5 pt,label=182: $v_3$] at (1,3) {} ;
		\node [fill=black, inner sep=2.5 pt,label=0: $v_4$] at (1.7,3) {} ;
		\node [fill=black, inner sep=1.5 pt] at (0,2) {} ;
		\node [fill=black, inner sep=1.5 pt] at (-1,2) {} ;
		\node [fill=black, inner sep=1.5 pt] at (1,2) {} ;
		\draw[line width=0.05cm, black] (0,1) -- (1,2)  {};
		\draw[line width=0.05cm, black] (0,1) -- (0,2)  {};
		\draw[line width=0.05cm, black] (0,1) -- (-1,2)  {};
		\draw[line width=0.05cm, black] (-1,2) -- (-1.8,3)  {};
		\draw[line width=0.05cm, black] (-1,2) -- (-0.6,3)  {};
		\draw[line width=0.05cm, black] (1,2) -- (1,3)  {};
		\draw[line width=0.05cm, black] (1,2) -- (1.7,3)  {};
		\draw [line width=0.04 cm, red] plot [smooth cycle] coordinates {(-1.8,3) (-2.2,5) (-1.3,5)};
		\draw [line width=0.04 cm, red] plot [smooth cycle] coordinates {(-0.6,3) (-0.9,4.5)  (-0.3,4.5)};
		\draw [line width=0.04 cm, red] plot [smooth cycle] coordinates {(1,3) (0.4,5.5) (1.4,5.5)};
		\draw [line width=0.04 cm, red] plot [smooth cycle] coordinates {(1.7,3) (1.7,5) (2.4,5)};
		\node[label=$T_1$] at (-1.83,3.7) {};
		\node[label=$T_2$] at (-0.6,3.65) {};
		\node[label=$T_3$] at (0.95,4.2) {};
		\node[label=$T_4$] at (1.95,4) {};
\end{tikzpicture}}

\usepackage{float}

\usepackage{color,soul}

\newcommand{\bb}{\begin{equation}}
	\newcommand{\ee}{\end{equation}}

\newcommand{\half}{\frac{1}{2}}

\newcommand{\cT}{\mathcal T}
\newcommand{\cB}{\mathcal B}
\newcommand{\cA}{\mathcal A}

\newcommand{\W}{\mathbb{W}}

\definecolor{dblue}{rgb}{.61,.61,1}
\definecolor{lightblue}{rgb}{.61,.61,1}
\definecolor{altblue}{rgb}{.61,.61,1}

\tikzset{cross/.style={cross out, draw=black, minimum size=2*(#1-\pgflinewidth), inner sep=0pt, outer sep=0pt},
	cross/.default={1pt}}

\newtheorem{thm}{Theorem}[section]
\newtheorem{prop}[thm]{Proposition}
\newtheorem{lem}[thm]{Lemma}
\newtheorem{cor}[thm]{Corollary}
\newtheorem*{thm*}{Theorem}
\theoremstyle{remark}

\theoremstyle{definition}

\def\IB{\relax\hbox{$\inbar\kern-.3em{\rm B}$}}
\def\IC{\relax\hbox{$\inbar\kern-.3em{\rm C}$}}
\def\ID{\relax\hbox{$\inbar\kern-.3em{\rm D}$}}
\def\IE{\relax\hbox{$\inbar\kern-.3em{\rm E}$}}
\def\IF{\relax\hbox{$\inbar\kern-.3em{\rm F}$}}
\def\IG{\relax\hbox{$\inbar\kern-.3em{\rm G}$}}
\def\IGa{\relax\hbox{${\rm I}\kern-.18em\mathcal{T}$}}
\def\IH{\relax{\rm I\kern-.18em H}}
\def\IK{\relax{\rm I\kern-.18em K}}
\def\IL{\relax{\rm I\kern-.18em L}}
\def\IP{\relax{\rm I\kern-.18em P}}
\def\IR{\relax{\rm I\kern-.18em R}}
\def\IZ{\relax{\rm Z\kern-.5em Z}}

\begin{document}
	\tikzset{
		position label/.style={
			below = 3pt,
			text height = 1.5ex,
			text depth = 1ex
		},
		brace/.style={
			decoration={brace, mirror},
			decorate
		}
	}

	\thispagestyle{empty}
	\quad
	
	\vspace{2cm}
	\begin{center}
		
		\textbf{\huge Trees with power-like height dependent weight}
		
		\vspace{1.5cm}
		
		{\Large Bergfinnur Durhuus ~~~~~~~~ Meltem Ünel}
		\vspace{0.5cm}
		
		{\it Department of Mathematical Sciences, Copenhagen University\\
			Universitetsparken 5, DK-2100 Copenhagen {\O}, Denmark}

		\vspace{0.5cm} {\sf durhuus@math.ku.dk, meltem@math.ku.dk}

		\vspace{1.5cm}
		
		{\large\textbf{Abstract}}\end{center} We consider planar rooted random trees whose distribution is even for fixed height $h$ and size $N$ and whose height dependence is given by a power function $h^\alpha$. Defining the total weight for such trees of fixed size to be $Z_N$, a detailed analysis of the analyticity properties of the corresponding generating function is provided. Based on this, we determine the asymptotic form of $Z_N$ and show that the local limit at large size is identical to the Uniform Infinite Planar Tree, independent of the exponent $\alpha$ of the height distribution function. 
	
	\vspace{0.5 cm}	
	

	\vspace{3cm}
	
	\newpage
	
	\tableofcontents	
	
	\section{Introduction}
	The study of random geometric objects has been actively pursued in recent years. In particular, deep results have been obtained for various models of random planar maps (or surfaces) concerning local limits and scaling limits. Thus, the Uniform Infinite Planar Triangulation was constructed in \cite{angel2003uniform} and the  Uniform Infinite Planar Quadrangulation appeared in \cite{chassaing2006local}, see also \cite{krikun2005local} and \cite{curien2012view,menard2010two}, while initial studies of the scaling limit of planar maps in \cite{chassaing2004random}  have been continued in numerous papers leading to proofs of existence of the limiting measure as well as establishing important properties of the limit, see e.g. \cite{marckert2006limit,le2007topological,le2012scaling}. While the subject is of natural interest in its own right within probability theory, important motivations and applications also arise from problems in theoretical physics, in particular statistical mechanics and quantum gravity, see e.g.  \cite{ambjorn1997quantum} and references therein.  
	
	A more classical topic in the same realm, and equally relevant from a physical point of view, is that of random trees, see \cite{drmota2009random} and references therein. Local limits and scaling limits have been constructed in this case as well and detailed information about the limiting measures can frequently be obtained using the highly developed theory of branching processes. 
	We shall in this paper concentrate on local limit results, in particular relating to the Uniform Infinite Planar Tree (UIPT) that can be obtained as a local (or weak) limit of uniformly distributed finite rooted planar trees of fixed size tending to infinity \cite{durhuus2003probabilistic}. Indeed, this limit can be viewed as a particular instance of a more general local limit result for critical Bienaym{\'e}-Galton-Watson (BGW)  trees conditioned on size, namely when the off-spring distribution is a geometric sequence, a result first established in \cite{kennedy1975galton} and \cite{aldous1998tree}. Earlier results on local limits of critical (and subcritical) BGW trees conditioned on height go back to H. Kesten \cite{kesten1986subdiffusive} and turn out to yield the same limit distribution supported on one-ended trees. More recent results allowing more general conditionings can be found in \cite{janson2012simply} and \cite{abraham2015introduction}. Properties of the limiting measures, e.g. relating to their Hausdorff dimension and spectral dimension, have been established in \cite{barlow2006random} and \cite{durhuus2007spectral}. 
	
	While the results just described depend heavily on the fact that the weights of individual trees are local, in the sense of being products of weights associated with the vertices of the tree, the case of non-local weight functions has been much less addressed in the literature. Thus, while the average height of various ensembles of planar trees has been of interest in many works, such as \cite{de1972average,meir1978altitude,flajolet1982average,chassaing2000height, richard2009onq},  it seems that properties of random planar trees with height-dependent weights have not been extensively explored, an exception being \cite{leckey2020height}, where a model of depth weighted random recursive trees is studied, with branching probability depending on vertex height.  
	
	In the present work we investigate a rather different case where the weight function $f(h)$ depends on the height $h$ of the tree and otherwise is constant for fixed size (see section 2 for a precise definition). It is easy to see that different limits can be obtained by a judicious choice of $f$. A detailed study of such cases will be given in \cite{meltem2021height}. Here, we consider the  case where $f$ is a power function, $f(h) = h^\alpha$, making use of transfer theorems from analytic combinatorics. The main result, stated in Theorem 4.2, is that the local limit in this case equals the UIPT, independently of the value of the exponent $\alpha$. 
	
	The paper is organised as follows. Section 2 contains a precise definition of the metric and associated Borel algebra on the space of rooted planar trees to be considered. Furthermore, the finite size distributions, whose limits we aim at calculating, are defined. In section 3, the analytic behaviour of the generating function for the total weights of trees as a function of the size $N$ is examined for fixed $\alpha$, based on well known results for the corresponding generating functions for trees of fixed height. Some technical details of the calculations involved are deferred to an appendix. Applying transfer theorems, the asymptotic behaviour of the coefficients at large $N$ is determined. These results are used in section 4 to prove existence of the weak limit and to identify it as the UIPT. Finally, section 5 contains some concluding remarks.

	\section{Preliminaries}\label{sec:2}
	
	In the following, ${\mathcal T}_N$ will denote the set of rooted planar trees of size $|T|=N$ with root vertex of degree $1$. Here $N\in{\mathbb N} :=\{1,2,3,\dots\}$ or $N=\infty$ (in which case $T$ is assumed to be locally finite) and we set 
	$$
	{\mathcal T} = \bigcup\limits_{N=1}^\infty {\mathcal T}_N\cup \cT_{\infty}\,,\quad \cT_{\rm fin} = \bigcup\limits_{N=1}^\infty {\mathcal T}_N\,.
	$$ 
	For $T\in\cT$ the notation $h(T)$ will be used for the height of $T$, i.e. the maximal length of a simple path in $T$ originating from the root. 
	For fixed $\alpha\in\mathbb R$ we consider the probability measures $\mu_N,\,N\in{\mathbb N},$ on $\cT$ given by
	\bb\label{probmeasure} 
	\mu_N(T) = Z_N^{-1}\cdot h(T)^{\alpha}\,,\quad T\in \cT_N,
	\ee
	where $Z_N$ is a normalisation factor, also called the \emph{finite size partition function}, given by 
	$$
	Z_N = \sum_{T\in\cT_N} h(T)^\alpha\,.
	$$
	Thus, $\mu_N$ is supported on $\cT_N$ and our goal is to study the weak limit of $\mu_N$ as $N\to\infty$ as a probability measure on $\cT$. Here $\cT$ is considered as a metric space whose metric $dist$ is defined as follows. Denoting the root of $T\in\cT$ by $v_0$, the ball $B_r(T)$ of radius $r$ in $T$ around $v_0$ is by definition the subgraph of $T$ spanned by vertices at distance at most $r$ from $v_0$, i.e.
	$$
	V(B_r(T))=\{v\in V(T)\mid d_T(v_0,v)\leq r\} \, ,
	$$ 
	where $d_T$ designates the graph distance on $T$ and $V(G)$ denotes the vertex set of any given graph $G$. For $T,T'\in \cT$ we then set 
	\bb\label{eq:metric}
	\mbox{\rm dist}(T,T') = \inf\{\frac{1}{r}\mid B_r(T)=B_r(T'),\, r\geq 1\}\,.
	\ee
	It is easily verified that $\mbox{\rm dist}$ is a metric on $\cT$. In fact, it is an ultrametric in the sense that for any triple $T,T',T''$ of trees in $\cT$ we have 
	$$
	\mbox{\rm dist}(T,T'') \leq \mbox{max}\{\mbox{\rm dist}(T,T'),\mbox{\rm dist}(T',T'')\}\,.
	$$
	The ball of radius $s>0$ around a tree $T$ will be denoted by $\cB_s(T)$. If $s=\frac 1r$, where $r\in\mathbb N$, it is seen that
	\bb\label{def:ball}
	\cB_{\frac 1r}(T) = \cB_{\frac 1r}(T_0) = \{T'\in\cT\mid B_r(T') = T_0\}\,,
	\ee
	where $T_0 = B_r(T)$. For additional properties of the metric space $(\cT,{\rm dist})$, including the fact that it is a separable and complete metric space, see e.g. \cite{durhuus2003probabilistic}.  
	
	By definition, a sequence of probability measures $(\nu_N)_{N\in\mathbb N}$ on $\cT$ is \emph{weakly convergent} to a probability measure $\nu$ on $\cT$	
	\bb\nonumber
	\int \limits_{\mathcal T} f d\nu_N \rightarrow \int \limits_{\mathcal T}  f d\nu
	\ee
	as $N\rightarrow \infty$ for all bounded continuous functions $f$ on $\mathcal{T}$. We refer to \cite{billingsley2013convergence} for a detailed account of weak convergence of probability measures. 
	
	As mentioned earlier, the main result of this paper is a proof that the weak limit of the sequence $(\mu_N)_{N\in\mathbb N}$ exists and equals the UIPT, which will be described in more detail in section 4.  A basic ingredient in the proof is the generating function $\mathbb W_\alpha$ for the $Z_N$, defined by
	\bb\label{eq:Walpha1}
	\W_\alpha (g) = \sum_{N=1} ^\infty \sum_{T\in\cT_N} h(T)^\alpha g^{|T|} = \sum_{N=1}^{\infty} Z_N g^N\,,
	\ee
	where $g$ is a real or complex variable. It is known, and will also be shown below, that the sum is convergent for $|g|<\frac 14$ and divergent for $|g|>\frac 14$. 
	It will be convenient for the following discussion to define
	\bb \nonumber \label{partfuncfinite}
	X_m(g) = \sum_{h(T)\leq m} g^{|T|}
	\ee
	for $m\geq 1$, such that 
	\bb\label{eq:Walpha2}
	\W_\alpha  = \sum_{m=1}^\infty m^\alpha (X_{m+1}-X_m)\,.
	\ee
	It is a known fact, and easy to verify, that $X_m$ fulfills the recursion relation
	$$
	X_{m+1}(g) = \frac{g}{1-X_m(g)}\,, m\geq 1, \quad \mbox{and}\quad X_1(g)=g\,.
	$$
	This equation can be rewritten in linear form and hence solved explicitly. The result is given by (see e.g. \cite{drmota2009random})
	\bb\label{eq:sol1}
	X_m(g) = 2g\frac{(1+\sqrt{1-4g})^m -(1-\sqrt{1-4g})^m}{(1+\sqrt{1-4g})^{m+1} - (1-\sqrt{1-4g})^{m+1}} \,,\quad m\geq 1\,.
	\ee
	Defining $X_0=0$ and 
	\bb\label{def:z}
	z=\sqrt{1-4g}\,,
	\ee
	this gives 
	\bb\label{eq:sol2}
	X_m(g) -X_{m-1}(g) = \frac{z^2}{\frac{1+z}{2} (\frac{1+z}{1-z})^m + \frac{1-z}{2} (\frac{1-z}{1+z})^m -1} \,,\quad m\geq 1\,.
	\ee
	Noting that for $z=x+iy\in \mathbb C$ we have
	$$
	\Big|\frac{1-z}{1+z}\Big|^2 = \frac{(1-x)^2+y^2}{(1+x)^2+y^2} \, ,
	$$
	it follows that $|\frac{1-z}{1+z}|<1$ if $x>0$, and hence, by \eqref{eq:sol1}, $X_m$ is an analytic function of $z$ away from the imaginary axis for all $m$. Note also that for fixed such $z$ the right-hand side of \eqref{eq:sol2} decays exponentially with $m$ and hence the series in \eqref{eq:Walpha2} converges. Taking $z$ to be real we get, in particular, that $\mathbb W_\alpha$ is well-defined for $0\leq g<\frac 14$. On the other hand, considering the denominator in \eqref{eq:sol1} we see that it is an odd polynomial in $z$ of degree $m+1$ if $m$ is odd and otherwise of degree $m$, and it is straightforward to see that its roots are given by 
	\bb \nonumber \label{def:zp}
	z_p = \pm i\tan\Big(\frac{p}{m+1}\pi\Big)\,,\quad p=0,1,\dots,\lfloor \tfrac{m}{2}\rfloor\,.
	\ee
	Hence, it follows from \eqref{eq:sol1} that these are precisely the (simple) poles of $X_m$ apart from $z=0$ corresponding to $p=0$ where both numerator and denominator have simple zeroes. Therefore, $X_m$ is analytic for $|z|<|z_1|$ and hence for $|g|<g_m$ where 
	\bb\label{def:gm}
	g_m= \frac 14\Big(1+\tan^2\frac{\pi}{m+1}\Big)\,,
	\ee
	with a singularity at $g=g_m$. Since $g_m\to \frac 14$ as $m\to \infty$ it follows from \eqref{eq:Walpha2} that the sum defining $\mathbb W_\alpha(g)$ is divergent for $g>\frac 14$, as claimed above.

	\section{Generating functions}
	
	\subsection{Singular behaviour}
	In order to establish the existence of the limit of $(\mu_N)$ we shall determine the asymptotic behaviour of $Z_N$ as $N\to\infty$. This will be done by analysing the singularity of $\mathbb W_\alpha$ at $g=\frac 14$ in more detail and using so-called transfer theorems \cite{flajolet2009analytic}.  When considering $W_\alpha$ as a function of $z$ we will use the notation $\W_\alpha(z)$ for $\W_\alpha(g)$. 
	
	Before stating the main result of this section, we note a few elementary facts that will be needed, formulated in the following three lemmas.
	
	\begin{lem} 
		For $|z| < 1$ and $m \geq 1 $ we have
		
		\begin{equation} \nonumber 
			\Big(\frac{1+z}{1-z} \Big) ^m = 1 + \sum_{k=1} ^{\infty} a_k ^m (2mz)^k\,,
		\end{equation}
		where the coefficients $a_k^m$ fulfill 
		\bb\label{eq:a_k1}
		0< a_k ^m < e^2 
		\ee
		for all $k$. Moreover,
		\bb \label{a_m}
		\abs{a_k ^m - \frac{1}{k!}} \leq \frac{e^{k+1}}{m}\quad\mbox{for $k\leq m$}.
		\ee

		\label{lem: lemma1}
	\end{lem}
	
	\begin{proof} First, we note that the term of order $k$ in the series expansion 
		\begin{eqnarray}
			\nonumber \Big( \frac{1+z}{1-z} \Big)^m &=& (1+ 2(z+z^2+...))^m\\
			\nonumber &=& 1 + \sum_{r=1} ^m \binom{m}{r} 2^r (z + z^2 +...)^r \\
			&=& 1+ \sum_{r=1} ^m \frac{(1-\frac{1}{m})\cdots(1-\frac{r-1}{m})}{r!} (2mz)^r (1+z+...)^r
			\label{eq: exp} \nonumber
		\end{eqnarray}
		is identical to the one in the polynomial
		\begin{eqnarray} \nonumber 
			P_k(z) &=& 1 + \sum_{r=1} ^m \frac{(1-\frac{1}{m})\cdots(1-\frac{r-1}{m})}{r!} (2mz)^r (1+z+...+z^{k-1})^r \\
			\nonumber &=& 1+ \sum_{r=1} ^m (1-\frac{1}{m}) \cdots (1-\frac{r-1}{m})(2mz)^r \sum_{r_0+...+r_{k-1}=r} \frac{1}{r_0 ! \cdots r_{k-1}!} z^{r_1} z^{2r_2}\cdots z^{(k-1)r_{k-1}}.
		\end{eqnarray}  
		From this expression we see that the coefficient of $(2mz)^k$ is given by	
		\begin{eqnarray}
			\label{am1} a_k^m &=& \sum_{r=1} ^m (1-\frac{1}{m})\cdots(1-\frac{r-1}{m}) (2m)^{r-k} \sum\limits_{\substack{r_0+...+r_{k-1}=r \\ r_1 + ...+(k-1)r_{k-1} = k-r}} \frac{1}{r_0! r_1 ! \cdots r_{k-1} !}\\
			\nonumber &<&  \sum_{r=1} ^m \frac{1}{r!} ( 1+ (2m)^{-1} +...+(2m)^{-(k-1)})^r\\\nonumber &<&\sum_{r=1}^m  \frac{1}{r!} \Big( \frac{2m}{2m-1} \Big)^r < e^2 \, ,
		\end{eqnarray}
		where the first inequality follows by replacing the prefactors by 1 and lifting the restriction $r_1+2r_2+...+(k-1)r_{k-1} = k-r$, while the second inequality results from estimating the finite sum in parenthesis by the corresponding infinite sum. This proves the first part of the lemma.
		
		To show \eqref{a_m}, we rewrite \eqref{am1} for $k\leq m$ in the form
		\begin{eqnarray}\nonumber
			a_k ^m &=&  \frac{1}{k!} (1- \frac{1}{m})\cdots(1-\frac{k-1}{m})  \\
			&+&  \sum_{r=1} ^{k-1} (1-\frac{1}{m})\cdots(1-\frac{r-1}{m}) (2m)^{r-k} \sum\limits_{\substack{r_0+...+r_{k-1}=r \\ r_1 + ...+(k-1)r_{k-1} = k-r}} \frac{1}{r_0! r_1 ! \cdots r_{k-1} !}\,. \nonumber
		\end{eqnarray}
		Here, the latter sum can be estimated as above by 
		\begin{equation}
			\nonumber\sum \limits_{r=1} ^{k-1} (\frac{1}{2m})^{k-r} \frac{1}{r!} \sum \limits_{r_0+..+r_{k-1} =r} \frac{r!}{r_0!\cdots r_{k-1}!} \;\leq \; \frac{1}{2m} \sum \limits_{r=1} ^{k-1} \frac{k^r}{r!} \;< \; \frac{e^k}{2m}\,,
		\end{equation}
		while expanding the first term yields an expression of the form $\frac{1}{k!} + \frac{c_k^m}{m}$ where the coefficient $c_k^m$ is easily seen to fulfill the bound
		$$
		|c_k^m|\leq \frac{2^k}{k!} < \frac{e^2}{2}  \,.
		$$  
		
		\noindent Combining the two estimates yields the claimed bound. \end{proof}
	
	The next lemma concerns properties of the denominator function
	\bb \nonumber  \label{def:Dm}
	D_m(z) = \frac{1+z}{2} \Big(\frac{1+z}{1-z}\Big)^m + \frac{1-z}{2} \Big(\frac{1-z}{1+z}\Big)^m -1\,,
	\ee
	appearing in \eqref{eq:sol2}, for small values of $z$. We use the notation $A_k$ for the Taylor coefficients of $\frac 12 t^2(\cosh t-1)^{-1}$, i.e.
	\bb \label{laurent}
	\frac{1}{\cosh(2t)-1} = \frac{1}{2t^2} \big[ 1+ \sum_{k=1} ^\infty A_k(2t)^{2k} \big]\,,\quad 0<|t|<\pi.
	\ee
	\begin{lem}
		a)\;\;	For $|z| <1$ and $m \geq 1$ it holds that
		\bb \nonumber  \label{eq:lem3a}
		D_m(z) = 2m(m+1) \cdot z^2 \cdot \Big( 1+ \sum_{k \geq 1} b^{m} _{2k} (2mz)^{2k} \Big)\,,
		\ee
		where the coefficients $b^m _{2k}$ fulfill
		\bb\label{eq:b_2k1}
		|b_{2k} ^m| \leq K\,,\quad k\geq 1 \, ,
		\ee
		as well as
		\bb \label{eq:b_2k2}
		\Big| b^m _{2k} - \frac{2}{(2(k+1))!}\Big|\, \leq\, L\cdot\frac{e^{2k}}{m}\,,\quad k\leq m \, ,
		\ee
		where $K$ and $L$ are positive constants. 
		\vspace{0.1 cm}
		
		b)\;\; For $|z| \leq \tan (\frac{\pi}{m+1})$ we have
		\bb \nonumber 
		\frac{z^2}{D_m(z)} = \frac{1}{2m(m+1)} \Big( 1+ \sum_{k \geq 1 } c_{2k}^m(2mz)^{2k} \Big)\,,
		\ee
		where  the coefficients $c_{2k}^m$ fulfill
		\bb\label{eq:lem3b}
		|c_{2k} ^m| \leq (K+1)^k\, , \quad  k\geq 1 \, ,
		\ee
		and there exist positive numbers $B_k, k=1,2,\dots$, independent of  m, such that 
		\bb \label{eq:c_2k2}
		\abs{ c_{2k} ^m - A_k} \leq  \frac{B_k}{m}\,,\quad k\leq m\,.
		\ee
		\label{lem: lemma2}
	\end{lem}
	
	\begin{proof}
		With notation as in Lemma \ref{lem: lemma1} we have
		
		\begin{eqnarray}
			\nonumber D_m(z) &=& \sum_{k \geq 1} \big(a_{2k}^m (2mz)^{2k} + z a_{2k-1} ^m (2mz)^{2k-1}\big)\\
			\nonumber &=& \sum_{k \geq 1} \big(a_{2k}^m + \frac{1}{2m}a_{2k-1}\big)(2mz)^{2k} \\
			\nonumber &=& 2m(m+1)z^2 \Big( 1 + \tfrac{1}{m+1}\sum_{k \geq 1} \big(2m a_{2(k+1)} ^m + a^m _{2k+1}\big) (2mz)^{2k} \Big)\,,
		\end{eqnarray}
		where it has been used that $a_1^m=1$ and $a_2^m=\frac 12$. Hence,
		\begin{align*}
			& b_{2k} ^m = \tfrac{1}{m+1}\big(2m a_{2(k+1)} ^m + a^m _{2k+1}\big) 
		\end{align*}
		and \eqref{eq:b_2k1} follows immediately from \eqref{eq:a_k1} with $K=2e^2$. Similarly, the bound \eqref{eq:b_2k2} follows easily from \eqref{eq:a_k1} and \eqref{a_m} with $L=2(e^3+1)$.
		This proves part a) of the lemma.
		
		For the second part we note that	$\frac{z^2}{D_m(z)}$ is a meromorphic function of $z$ which is analytic for $|z|< \tan (\frac{\pi}{m+1})$ as shown in section \ref{sec:2}. The power series for $z^2/D_m(z)$ in this disc is obtained by inverting that of $\frac{D_m(z)}{z^2}$, i.e. the coefficients $c^m_{2k}$ are determined by
		
		\begin{equation}\label{eq:c_2k3}
			c^m _{2k} = - \sum_{l=0} ^{k-1} c_{2l} ^m\cdot b^m _{2(k-l)}\,,
		\end{equation}
		where $c_0^m \equiv 1$ . Using the the bound \eqref{eq:b_2k1} just proven, we obtain
		
		\bb\nonumber
		|c_{2k} ^m| \leq K \sum_{l=0} ^{k-1} |c^m _{2l}|\,,
		\ee
		which implies \eqref{eq:lem3b} by a simple induction argument. 
		
		Using that  
		\bb \nonumber \label{Ak}
		A_k = - \sum _{l=0} ^{k-1} \frac{2A_l}{(2(k-l+1))!} 
		\ee
		where $A_0 \equiv 1$, it follows from \eqref{eq:c_2k3} that
		\bb \nonumber
		c_{2k} ^m - A_k = \sum_{l=0} ^{k-1}  \Big(\frac{2A_l}{(2(k-l+1))!} - c_{2l} ^m b^m _{2(k-l)}\Big)\,.
		\ee
		Using the bounds of part a), this gives 
		\bb \nonumber
		\mid c_{2k} ^m - A_k\mid  \leq \sum_{l=0} ^{k-1}  \Big(|A_l| L\frac{e^{2(k-l)}}{m} + K|c_{2l}^m-A_l|\Big)\,.
		\ee
		Since $c_0=A_0=1$, the inequality \eqref{eq:c_2k2} now follows by induction with $B_k$ defined recursively as
		$B_0=0$ and 
		$$
		B_{k+1} = \sum_{l=0}^{k} (L|A_l|e^{2(k-l)} + KB_l)\,, \quad k\geq 0.
		$$
		This concludes the proof.	\end{proof}
	
	For fixed $a>0$ we denote in the remainder of this paper by $V_a$ the wedge in the right half-plane given by
	$$
	V_a = \{z=x+iy \in\mathbb C\mid x>0,\, |y|<ax\}\,.
	$$  
	
	\begin{lem} Let $a>0$ and $0<\epsilon<1$ be given. Then there exist positive constants $K_0, K_1, m_0$ depending on $a$, and $\delta_0$ depending on $a$ and $\epsilon$, such that the following statements hold.
		\begin{enumerate}[label=\alph*)]
			\item For $z=x+iy \in V_a $,
			\bb
			2 \sinh ^2(mx) \leq  |\cosh (2mz)-1| \leq K_0 \sinh ^2 (mx)\,.
			\label{eq:lem4a}
			\ee
			
			\item For $z=x+iy \in V_a$, $|z| \leq  \delta_0$, and $ m \geq m_0$,
			\bb
			\sinh ^2((1-\epsilon)mx) \leq |D_m(z)| \leq K_1 \sinh ^2 ((1+\epsilon)mx)\,.
			\label{eq:lem4b}
			\ee
		\end{enumerate}		
		\label{lem: bounds}
	\end{lem}
	
	\begin{proof} 
		~~ \newline
		$a)$ \;\; For any $z\in\mathbb C$ we have
		
		\begin{equation}\nonumber |\cosh (2z)-1| 
			= 2\big( \sinh^2 x + \sin ^2 y\big) \geq 2 \sinh^2  x\,.
		\end{equation}
		On the other hand, for $z \in V_a$ we have
		
		\bb\nonumber
		\sin^2 (my) \leq (my)^2 \leq a^2 (mx)^2 \leq a^2 \sinh^2 (mx)\,.
		\ee
		Combining these two estimates yields \eqref{eq:lem4a} with $K_0=2(1+a^2)$.
		
		\noindent $b)$\;\; First, note that
		
		\bb \nonumber 
		\Big( \frac{1+z}{1-z} \Big) ^m = e^{2mz(1+ O(z))}\,,
		\label{def:O}
		\ee
		where $O(z)$ is analytic and fulfills 
		\bb \nonumber  \label{eq:Oineq}
		|O(z)|\leq c|z|
		\ee
		for  $z$ in a suitably small disc around $0$, where $c$ is some constant independent of $m$. Therefore, choosing  $\delta_0$ small enough, it follows that $2mz(1+ O(z)) \in V_{2a}$ if $|z|<\delta_0\,, z\in V_a$. Introducing the shorthand 
		\bb\label{def:zprime}		
		z^\prime=x^\prime+iy^\prime= z(1+ O(z))\,,
		\ee
		we get
		\begin{eqnarray}
			\nonumber 
			|D_m(z)| &=& |\frac{1+z}{2} e^{2mz^\prime} + \frac{1-z}{2} e^{-2mz^\prime} -1| \\
			\nonumber &=& |\cosh(2mz^\prime) + z \sinh(2mz^\prime)-1| \\
			\nonumber &\geq& |\cosh(2mz^\prime)-1| -|z| |\sinh(2mz^\prime)| \\
			\nonumber &\geq& 2\sinh^2 (mx^\prime) - |z|(\sinh^2 (2mx^\prime) + \sin^2 (2my^\prime))^\half \\ 
			&=& 2 \sinh^2 (m x^\prime) \Big[ 1- \frac{|z|}{2 \sinh(mx^\prime)} \Big( \frac{\sinh^2 (2mx^\prime)}{\sinh^2 (mx^\prime)} + \frac{\sin^2 (m y^\prime)}{\sinh^2(m x^\prime)} \Big)^\half \Big]\,.\label{lowDm}
		\end{eqnarray}
		Consider first the last expression in the case where $mx^\prime \geq 1$. In the second term inside square brackets, the factor multiplying $|z|$ is numerically bounded by a constant depending only on $a$. Hence, choosing $\delta_0$ sufficiently small, the term in square brackets is bounded from below, say by $\half$, for $|z|<\delta_0$. 
		
		On the other hand, if $mx^\prime \leq 1$ we observe that the factor inside  round brackets is bounded while the prefactor $\frac{|z|}{\sinh(mx^\prime)}$ can be estimated as follows. First, choosing $\delta_0$ sufficiently small such that $|O(z)|\leq \frac 12$ for $|z|<\delta_0$ we get     
		$$
		|z| \leq \frac{|z^\prime|}{1-|O(z)|}\leq 2|z^\prime|\leq 2(1+2a)x^\prime
		$$
		and hence 
		$$
		\frac{|z|}{\sinh(mx^\prime)} \leq \frac{2(1+2a)}{m}\frac{mx^\prime}{\sinh(mx^\prime)} \leq \frac{2(1+2a)}{m}\,.
		$$
		Then, choosing $m_0$ sufficiently large, the term in square brackets in \eqref{lowDm} is bounded from below by $\frac 12$ if $m\geq m_0$.   
		Thus, we have shown that for $m_0$ sufficiently large and $\delta_0>0$ small enough it holds that
		\bb
		|D_m(z)| \geq \sinh^2(mx^\prime)\quad\mbox{for $z\in V_a,\, |z|<\delta_0,\, m>m_0\,.$}
		\label{eq: ineqk1}
		\ee
		Given $\epsilon$, choose $\delta_0$ small enough such that $|zO(z)|< \frac{|z|}{1+a}\epsilon$ for $|z|<\delta_0$. For such $z$ in $V_a$ we then have 
		\bb\nonumber
		x^\prime = x + \operatorname{Re}(zO(z)) = x(1+ \operatorname{Re}(\frac{z}{x}O(z)))\,,
		\ee
		where 
		$$
		|\operatorname{Re}(\frac{z}{x}O(z))| \leq \frac{|z|}{(1+a)x}\epsilon <\epsilon
		$$		
		and hence
		\bb \nonumber
		(1- \epsilon)x \leq x^\prime \leq (1+ \epsilon)x\,,
		\ee
		from which the lower bound in \eqref{eq:lem4b} follows in view of \eqref{eq: ineqk1}.
		
		By a slight modification of the arguments above, the upper bound follows similarly.
	\end{proof}

	\noindent We are now ready to prove the first main result of this section.
	
	\begin{thm}
		Assume the exponent $\alpha$ fulfills $-(2n+1) < \alpha < -(2n-1), \, n \in \mathbb{N}_0$,\footnote{We use ${\mathbb N}_0$ to denote the non-negative integers.} and let $a>0$. Then $\W_\alpha$ is analytic in the right half-plane ${\mathbb C}_+ =\{z\in {\mathbb C}\mid \operatorname{Re} z > 0\}$ and there exists a polynomial $\W^{(n)} _\alpha(z)$ of degree $2n$, a real constant $c_\alpha$, and $\Delta>0$ such that
		\bb 
		\W_\alpha (z) = \W_\alpha ^{(n)} (z) + c_\alpha z^{1-\alpha} + O(|z|^{1-\alpha+ \Delta})
		\label{eq: partfunc}
		\ee
		for $z$ small in $V_a$.
		\label{thm: main1}
	\end{thm}
	
	\begin{proof}
		We claim that the polynomial $\W_\alpha ^{(n)}$ is given by replacing the summand in the definition \eqref{eq:Walpha2} of $\W_\alpha$ by its Taylor polynomial of degree $2n$, i.e. 
		\bb \label{Wn}
		\W_\alpha ^{(n)} (z) = \sum_{m=1} ^{\infty} \frac{m^\alpha}{2m(m+1)} \big( 1+ \sum \limits_{k=1} ^n c_{2k} ^m (2mz)^{2k} \big)\,,
		\ee
		and that $c_\alpha$ is given by
		\bb \nonumber 
		c_\alpha = \int_{0} ^{\infty}  \Big( \frac{t^{\alpha}}{\cosh (2t)-1} - t^\alpha L_n(2t) \Big) dt\,, 
		\ee
		where $L_n(t)$ stands for the Laurent polynomial of order $2(n-1)$ for $\frac{1}{\cosh t-1}$, i.e. 
		\bb\label{eq:Ln}
		L_n(t) = \frac{2}{t^2}\Big[1+\sum_{k=1}^n A_kt^{2k}\Big]\,.
		\ee
		Setting
		\bb \label{eq:Tnm}
		T_n^m(z) = \frac{1}{2m(m+1)} \Big( 1+ \sum \limits_{k=1} ^n c_{2k} ^m (2mz)^{2k} \Big)\,,
		\ee
		we start by rewriting 
		\begin{align}
			\nonumber &\W_\alpha (z) - \W_\alpha ^{(n)}(z) - c_\alpha z^{1-\alpha}  \\
			\nonumber &= \sum \limits_{m=1} ^\infty \Big[ \frac{z^2 m^\alpha}{D_m(z)} - m^\alpha T_n^m(z) \Big] - z^{1-\alpha}\int \limits_0 ^\infty \Big( \frac{t^\alpha}{\cosh(2t)-1} - t^\alpha L_{n}(2t) \Big)dt \\
			&= \sum \limits_{m=1} ^\infty \Big[ \frac{z^2 m^\alpha}{D_m(z)} - m^\alpha T_n^m(z) \Big] - \int \limits_0 ^\infty \Big( \frac{z^2t^\alpha}{\cosh(2tz)-1} - z^{2} t^\alpha L_n(2tz)\Big) dt\,,
			\label{4piece}
		\end{align}
		where in the last line we converted the integral over the positive real axis to a line integral along the half line $\ell_z: t \rightarrow tz$ inside the wedge $V_a$, by using Cauchy's theorem and the fact that $\alpha+2k-2<-1$ for $k \in \{0,1,..,n\}$ implies  that the integrand decays fast enough and uniformly to $0$ at infinity inside $V_a$, such that the integral along the circular part of the contour at infinity vanishes.
		
		It is convenient to split the sum and integral in the last expression into three regions given by
		\begin{flalign}
			\nonumber & 1) ~ m|z|\leq r ~, ~ t|z|\leq r &&\\
			& \label{split} 2) ~ r< m|z|\leq r^{-1} ~,~ r< t|z|\leq r^{-1} &&\\
			\nonumber & 3) ~ r^{-1} < m|z| ~,~ r^{-1} < t|z|\,, &&
		\end{flalign} 
		
		\noindent respectively\footnote{More precisely, $t$ is restricted according to $0<t\leq \big\lfloor\frac{r}{|z|}\big\rfloor$, $\big\lfloor\frac{r}{|z|}\big\rfloor<t \leq \big\lfloor\frac{1}{r|z|}\big\rfloor$ and $\big\lfloor\frac{1}{r|z|}\big\rfloor<t$, respectively, where $\lfloor x \rfloor$ denotes the integer part of $x\in\mathbb R$. For the sake of simplicity, we shall in the following apply the notation in \eqref{split} and omit writing explicitly the relevant integer parts.}. Here, $r$ will ultimately be chosen as a suitable function of $|z|$ tending to $0$ as $z$ tends to $0$, but for the moment we merely assume $0<r<1$. Denoting the corresponding contributions to the series and integral in \eqref{4piece} by $S_1, S_2, S_3$ and $I_1, I_2, I_3$, respectively, we have  
		$$
		\W_\alpha(z) -\W_\alpha ^{(n)} (z)  -c_\alpha z^{1-\alpha}  = S_1 + S_2 +S_3 - (I_1 + I_2 +I_3)\,,
		$$
		and we shall proceed by successively estimating the numerical values of $S_1, I_1$ first, then $S_2 - I_2 $, and finally of $S_3, I_3$. In the remainder of this proof, "cst" will denote a generic positive constant independent of $z$ and $r$ within the stated ranges, and likewise $O(w)$ will denote a generic function of order $w$ for $|w|$ small, i.e. $|O(w)|\leq \text{cst}\cdot |w|$.
		
		\bigskip
		
		\noindent \underline{$S_1$ and $I_1$:} \; It follows from Lemma \ref{lem: lemma2} that
		\bb\label{eq:r1}
		\frac{z^2}{D_m(z)} = T_n^m(z) + (mz)^{2n}O(z^2)
		\ee
		provided $r$ fulfills
		\bb\label{eq:cond1}
		4r^2(K+1) < 1\,,
		\ee
		which we henceforth assume to hold. It follows that
		
		\begin{align}
			\nonumber
			|S_1| ~& \leq ~\sum \limits_{m|z|\leq r} \Big|\frac{z^2m^\alpha}{D_m(z)}- m^\alpha T_n^m(z)\Big| \\
			&= \sum \limits_{m|z|\leq r} (m|z|)^{2n+\alpha} O(|z|^{2-\alpha}) \leq \text{cst} \cdot |z|^{1-\alpha}r^{1+\alpha+2n}.
			\label{eq:estS1}
		\end{align}
		
		\noindent For $I_1$, on the other hand, we get 
		\begin{align}
			\nonumber |I_1| \leq \int \limits_0 ^{\frac{r}{|z|}} \big| \frac{z^2t^\alpha}{\cosh(2tz)-1}-z^2t^{\alpha}L_n(2tz) \big| dt 
			&\leq  \int \limits_0 ^{\frac{r}{|z|}} t^{\alpha}|z|^2 O(|tz|^{2n})dt \\
			& \leq \text{cst} \cdot |z|^{1-\alpha}r^{1+\alpha+2n}\,.
			\label{eq:estI1}
		\end{align}
		
		\bigskip
		
		\noindent \underline{$S_2-I_2$:}\; We further decompose this expression as 
		$$
		S_2-I_2 = A+ B+ C\,,
		$$
		where 
		\begin{align}
			\label{A} A & = \sum_{r<m|z| \leq \frac{1}{r}} \Big( \frac{z^2 m^\alpha}{D_m(z)} - \frac{z^2 m^\alpha}{\cosh(2mz)-1} \Big)\\ \label{B} B & =  z^{2} \Big[ \sum_{r<m|z| \leq \frac{1}{r}}  \frac{ m^\alpha}{\cosh(2mz)-1} - \int \limits_{\frac{r}{|z|}} ^{\frac{1}{r|z|}} \frac{t^\alpha dt}{\cosh(2tz)-1} \Big]\\  \label{C}  C &= \int \limits_{\frac{r}{|z|}} ^{\frac{1}{r|z|}} z^{2} t^\alpha L_n(2tz) dt - \sum \limits_{r<m|z|\leq \frac{1}{r}} m^\alpha T_n^m(z)\,,  
		\end{align}
		
		\noindent and proceed to establish bounds on $A, B$ and $C$. In order to deal with $A$, we note that with the notation $z^\prime$ introduced earlier in \eqref{def:zprime} we have
		
		\begin{equation}\label{eq:diffA}
			\frac{1}{\cosh(2mz) -1}-\frac{1}{D_m(z)} = \frac{\cosh(2mz^\prime)- \cosh (2mz) + z \sinh(2mz^\prime)}{D_m(z)(\cosh(2mz)-1)}\,.
		\end{equation}
		By Lemma \ref{lem: bounds} with, say, $\epsilon =\frac 12$ the denominator fulfills
		
		\bb \label{denom}
		|D_m(z)(\cosh(2mz)-1)| \geq \text{cst} \cdot (mx)^4\quad\mbox{if $r <m|z| < 1$}
		\ee
		and 
		\bb\label{eq:cond2}
		|z| <  \min\{\delta_0, r/m_0\}\,.
		\ee 
		To estimate the numerator, we note the inequalities
		\begin{equation} \label{eq:ineq}
			|e^{2mz^\prime}- e^{2mz}| \leq  \text{cst} \cdot e^{2mx} m|z|^2 \quad\text{and}\quad 
			|e^{-2mz^\prime}-e^{-2mz}| \leq \text{cst} \cdot m|z|^2\,,
		\end{equation}
		valid if $m|z|^2$ is bounded, which is the case if \eqref{eq:cond2} holds and $m|z|<\frac{1}{r}$. For $m|z|\leq 1$ we hence get 
		\begin{equation} \nonumber
			\abs{\cosh(2mz^\prime)-\cosh(2mz)}\leq \text{cst} \cdot m|z|^2\,,\quad |z \sinh(2mz^\prime)| \leq \text{cst} \cdot m|z|^2 \,,
		\end{equation}
		and together with \eqref{eq:diffA} and \eqref{denom} this implies 
		
		\bb \label{A1}
		\Big|\frac{1}{\cosh(2mz)-1} - \frac{1}{D_m(z)}\Big| \leq \text{cst} \cdot \frac{|z|}{r^3}\quad\mbox{for $r<m|z|<1$} \, ,
		\ee
		provided  \eqref{eq:cond2} holds, which we assume in the following.
		
		On the other hand, in the range $1<m|z| < r^{-1}$ the inequalities \eqref{eq:ineq} imply
		
		\begin{equation}\nonumber
			\nonumber |e^{2mz^\prime} - e^{2mz}| \leq \text{cst} \cdot e^{2mx} \frac{|z|}{r}\quad \text{and}\quad 
			|e^{-2mz^\prime}-e^{-2mz}| \leq \text{cst} \cdot \frac{|z|}{r}\,,
		\end{equation}
		which in turn give
		\bb \nonumber
		\abs{\cosh(2mz^\prime)-\cosh(2mz)}\leq \text{cst} \cdot e^{2mx}\frac{|z|}{r}\,,\quad |z| |\sinh(2mz^\prime)| \leq \text{cst} \cdot e^{2mx} \frac{|z|}{r}\,.
		\ee
		Combining this with the bound 
		\bb \nonumber
		|D_m(z) (\cosh(2mz)-1)|\geq \text{cst} \cdot \sinh^2(mx) \sinh^2(\tfrac{1}{2} mx),
		\ee
		which follows from Lemma \ref{lem: bounds} a) and b) with $\epsilon=\frac 12$, we get
		
		\bb
		\Big|\frac{1}{\cosh(2mz)-1} - \frac{1}{D_m(z)}\Big| \leq \text{cst} \cdot \frac{|z|}{r} \, .
		\label{A2}
		\ee
		
		\noindent Thus, using \eqref{A1} and \eqref{A2}, we obtain the bound
		
		\bb\label{eq:estA}
		|A| \leq  \sum_{r < |mz| \leq r^{-1}} \Big|\frac{z^2m^\alpha}{\cosh(2mz)-1} - \frac{z^2m^\alpha}{D_m(z)}\Big| \leq \text{cst} \cdot \frac{|z|^{2-\alpha}}{r^{4+|\alpha|}}\,,
		\ee
		valid for any (fixed) value of $\alpha\in\mathbb R$.
		\smallskip
		
		To derive an upper bound for $|B|$ it is useful to rewrite 
		\begin{align}
			\label{Bintegrand} - B &= \sum_{r < m|z|\leq r^{-1}}z^2\int_m^{m+1} dt \int_m ^t\frac{d}{ds} \frac{s^\alpha}{\cosh(2sz)-1}ds\\
			&= \sum_{r < m|z|\leq r^{-1}}z^2\int_m^{m+1} dt \int_m ^t\Big(\frac{\alpha s^{\alpha-1}}{\cosh(2sz)-1} - \frac {2z
				s^\alpha \sinh(2sz)}{(\cosh(2sz)-1)^2}\Big) ds\,.\nonumber
		\end{align}
		It is straight-forward to estimate the integrand using \eqref{eq:lem4a} and recalling that $\alpha<1$. One finds 
		\begin{equation} \nonumber
			\Big|\frac{d}{ds}\frac{s^\alpha}{\cosh(2sz)-1}\Big| \leq
			\begin{cases}
				\text{cst} \cdot \frac{|z|^{1-\alpha}}{r^{3-\alpha}}, & \text{if} \ r<|sz| \leq 1 \\ ~~ \\
				\text{cst} \cdot \frac{|z|^{1-\alpha}}{r}, & \text{if} \ 1<|sz|\leq r^{-1}\,,
			\end{cases}
		\end{equation}
		which yields the bound
		\begin{equation} \label{eq:estB}
			\abs{B}   \leq |z|^2(r|z|)^{-1}\Big[ \text{cst} \cdot \frac{\abs{z}^{1-\alpha}}{r^{3+|\alpha|}} \Big] 
			\leq \text{cst} \cdot \frac{|z|^{2-\alpha}}{r^{4+ |\alpha|}}\,.
		\end{equation}

		In order to estimate $|C|$, we first use the mean value theorem to write
		
		\begin{align}
			\nonumber |C| &\leq \sum_{r<|mz|\leq r^{-1}} \abs{\frac{m^\alpha}{2m(m+1)}\sum_{k=0}^n c^m_{2k} (2mz)^{2k} 
				- \int_m^{m+1}\frac{t^{\alpha-2}}{2}\sum_{k=0}^n A_k(2tz)^{2k}dt }\\
			&\leq \sum_{k=0}^n \sum_{r<|mz|\leq r^{-1}} \abs{A_k(2z)^{2k}\Big[\frac{m^{\alpha-1+2k}}{2(m+1)}-\frac{t_m ^{\alpha-2+2k}}{2}\Big]+ \frac{m^{\alpha-1}}{2(m+1)}(c^m_{2k}-A_k)(2mz)^{2k}}\,,\label{eq:estCprime} 
		\end{align}
		where $t_m \in [m,m+1]$. Using 
		$$
		\Big|\frac{m^{\alpha-1+2k}}{m+1}-t_m ^{\alpha-2+2k}\Big| \leq (3-\alpha-2k) m^{\alpha-3+2k}
		$$
		and recalling  \eqref{eq:c_2k2}, which holds for all $m$ in the summation range provided  $z$ fulfills the condition 
		\bb \label{cond3}
		|z|<\frac{r}{n}\,,
		\ee
		we see that \eqref{eq:estCprime} implies 
		\bb\label{eq:estC}
		|C|\leq \text{cst} \cdot \sum \limits_{k=0} ^n  \frac{|z|^{2-\alpha}}{r^{2-\alpha-2k}} \leq \text{cst} \cdot  \frac{|z|^{2-\alpha}}{r^{2-\alpha}} 
		\ee
		for such values of  $z$.		
		\smallskip
		
		From the estimates \eqref{eq:estA}, \eqref{eq:estB} and \eqref{eq:estC} we finally get 
		
		\begin{align}
			|S_2-I_2| \leq |A|+|B|+|C|\leq \text{cst} \cdot  \frac{|z|^{2-\alpha}}{r^{4+|\alpha|}} 
			\label{eq:est2}
		\end{align}
		if \eqref{eq:cond2} and \eqref{cond3} hold.
		
		\bigskip
		
		\noindent \underline{$S_3$ and $I_3$}:\; Using Lemma \ref{lem: bounds} a) we have
		\begin{equation}\label{S31}
			\Big| \int_{\frac{1}{r|z|}} ^\infty \frac{z^2 t^\alpha}{\cosh(2tz)-1} dt\Big| \leq  \text{cst} \int \limits_{\frac{1}{r|z|}} ^\infty \frac{|z|^2 t^\alpha}{\sinh ^2 (xt)} dt 
			\leq  \text{cst} \cdot \exp{\big(-\frac{\text{cst}}{r}\big)}|z|^{1- \alpha}\,,
		\end{equation}
		and similarly Lemma \ref{lem: bounds} b) yields
		\begin{equation}\label{I31}
			\Big|\sum_{m|z|>r^{-1}} \frac{z^2 m^\alpha}{D_m(z)}\,\Big|\, 
			\leq\, \text{cst} \cdot \exp{\big(-\frac{\text{cst}}{r}\big)} |z|^{1-\alpha}\,.
		\end{equation}
		Furthermore, using \eqref{eq:lem3b} we have
		\begin{align} \nonumber
			\Big|\sum \limits_{m|z|>r^{-1}} m^\alpha T_n^m(z)\,\Big| &\,\leq\, \sum \limits_{m|z|>r^{-1}} \sum_{k=0} ^n \frac {m^\alpha}{2m(m+1)}  (K+1)^k(2m|z|)^{2k} \\ &\leq \text{cst} \cdot |z|^{1-\alpha} r^{1-\alpha-2n}\,.  \label{Wn2}
		\end{align}
		Taking into account also the obvious bound
		\begin{align} \label{Tn}
			\Big|\,\int _{\frac{1}{r|z|}}^\infty z^{2} t^\alpha  L_n(2tz)dt\,\Big| \,\leq\, \text{cst} \cdot |z|^{1-\alpha} r^{1-\alpha-2n}\,,
		\end{align}
		it follows that
		\begin{equation} 
			|S_3|+|I_3| \leq \text{cst} \cdot |z|^{1-\alpha}r^{1- \alpha -2n} \label{eq:est3}\,.
		\end{equation}
		
		\smallskip
		
		Combining the estimates \eqref{eq:estS1}, \eqref{eq:estI1}, \eqref{eq:est2} and \eqref{eq:est3}, we conclude that
		
		\bb \label{eq:estall}
		\big|\W_\alpha(z) - \W_\alpha ^{(n)} - c_\alpha z^{1-\alpha}\big|\, \leq \,\text{cst} \cdot \Big(\frac{|z|^{2-\alpha}}{r^{4+|\alpha|}} + |z|^{1-\alpha} r^{1-\alpha-2n} \Big) \,,
		\ee
		provided $z$, $r$ and $n$ fulfill conditions \eqref{eq:cond1}, \eqref{eq:cond2} and \eqref{cond3}.  
		Choosing $r=|z|^\beta$ where $0<\beta<1$, these conditions are evidently satisfied for $|z|$ small enough. Noting that $1-\alpha-2n >0$, it follows from \eqref{eq:estall} that if $\beta$ is chosen such that    
		$$
		\beta(4+|\alpha|)<1\,,
		$$
		then \eqref{eq: partfunc} holds with $\Delta = \min\{1-\beta(4+|\alpha|),\beta(1-\alpha-2n)\}$.
		
		This concludes the proof of Theorem \ref{thm: main1}. 
	\end{proof}
	
	The next theorem is similar to the previous one and covers the case $\alpha >1 $.
	
	\begin{thm}
		
		\noindent Assume  $\alpha > 1$ and let $a>0$. Then $\W_\alpha$ is analytic in the right half-plane ${\mathbb C}_+$ and there exists  $\Delta>0$ such that
		\bb
		\W_\alpha (z) = c_\alpha z^{1-\alpha} + O(|z|^{1-\alpha+ \Delta})
		\label{eq: partfunc2}
		\ee
		for $z\in V_a$ small, where 
		\bb \nonumber 
		c_\alpha = \int_{0} ^{\infty} \frac{t^{\alpha}}{\cosh (2t)-1} dt\,.
		\ee
		\label{thm: main2}
	\end{thm}
	
	\begin{proof}
		Applying Cauchy's theorem as in \eqref{4piece} gives
		\bb \label{2piece}
		\W_\alpha (z) - c_\alpha z^{1-\alpha} = \sum \limits_{m=1} ^{\infty} \frac{z^2m^\alpha}{D_m(z)} - \int \limits_0 ^\infty \frac{z^2t^\alpha}{\cosh (2tz)-1}dt \, .
		\ee
		
		\noindent As previously, we split the sum and integral in \eqref{2piece} into three regions as in \eqref{split} and denote the corresponding contributions by $S_1,S_2,S_3$ and $I_1,I_2,I_3$, respectively. The relevant estimates can then be obtained by suitably modifying the arguments of the previous proof as follows. 
		\bigskip
		
		\noindent \underline{$S_1$ and $I_1$:} \; 
		Assuming $r$ fulfills \eqref{eq:cond1}, it follows from \eqref{eq:r1} that
		\begin{align} \label{S1I1i}
			\abs{S_1} \leq \; \sum \limits_{m|z|\leq r}\Big|\frac{z^2m^\alpha }{D_m(z)}\Big|\, \leq \sum \limits_{m|z|\leq r} \Big(\frac{m^\alpha}{2m(m+1)} + m^{\alpha-2}O((m|z|)^2)\Big) \leq \text{cst} \cdot |z|^{1-\alpha}r^{\alpha-1} \, .
		\end{align}
		Similarly, expanding the integrand of $I_1$ as in \eqref{laurent}, we get
		\begin{align} \label{S1I1ii}
			\abs{I_1} \leq \; \int \limits_0 ^{\frac{r}{|z|}} \frac{t^{\alpha-2}}{2}(1 + O((t|z|)^2) dt \, \leq \, \text{cst} \cdot |z|^{1-\alpha}r^{\alpha-1} \, .
		\end{align}
		
		\noindent \underline{$S_2-I_2$:} \; We further decompose this expression as
		\bb \nonumber
		S_2-I_2 = A+B \, ,
		\ee
		where 
		$A$ and $B$ are defined as in \eqref{A} and \eqref{B}.
		Since the estimates \eqref{A1} and \eqref{A2} do not depend on $\alpha$, the same bound as in \eqref{eq:estA},
		\bb
		|A| \leq \text{cst} \cdot \frac{\abs{z}^{2-\alpha}}{r^{4+\alpha}}\,, \label{Apart2}
		\ee
		holds in this case.
		
		Concerning $B$, we recall \eqref{Bintegrand} and estimate the integrand therein for $\alpha>1$ using \eqref{eq:lem4a}, which gives 
		\begin{equation} \nonumber
			\Big|\frac{d}{ds}\frac{s^\alpha}{\cosh(2sz)-1}\Big| \leq
			\begin{cases}
				\text{cst} \cdot \frac{|z|^{1-\alpha}}{r^{2}}, & \text{if} \ r<|sz| \leq 1 \\ ~~ \\
				\text{cst} \cdot \frac{|z|^{1-\alpha}}{r^\alpha}, & \text{if} \ 1<|sz|\leq r^{-1}\,.
			\end{cases}
		\end{equation}
		\noindent This implies the bound  
		\bb
		|B| \leq z^2 (r|z|)^{-1} [\text{cst} \cdot \frac{|z|^{1-\alpha}}{r^{\alpha}}] \leq \text{cst} \cdot \frac{|z|^{2-\alpha}}{r^{1+\alpha}}\,. \label{Bpart2}
		\ee
		\noindent From \eqref{Apart2} and \eqref{Bpart2} we get
		
		\bb 
		|S_2-I_2| \leq |A| + |B| \leq \text{cst} \cdot \frac{|z|^{2-\alpha}}{r^{4+\alpha}}
		\label{S2I2} \, .
		\ee
		
		\noindent \underline{$S_3-I_3$:} \; The estimates \eqref{S31} and \eqref{I31} are still valid for $\alpha>1$, hence we get
		\bb \label{S3I3}
		|S_3-I_3| \leq \text{cst} \cdot 
		\exp(- \frac{\text{cst}}{r})|z|^{1-\alpha}\,.
		\ee
		
		Collecting the estimates \eqref{S1I1i}, \eqref{S1I1ii}, \eqref{S2I2} and \eqref{S3I3}, we conclude that
		\bb \label{estall2}
		\abs{\W_\alpha(z) - c_\alpha z^{1-\alpha}} \leq \text{cst} \cdot \Big( \frac{|z|^{2-\alpha}}{r^{4+\alpha}} + |z|^{1-\alpha} r^{\alpha-1}\Big)\,,
		\ee
		
		\noindent provided $z$ and $r$ fulfill \eqref{eq:cond1} and \eqref{eq:cond2}.  
		Choosing $r=|z|^\beta$ as previously, where
		$ \beta(4+\alpha)<1, $
		it follows from \eqref{estall2} that \eqref{eq: partfunc2} holds with $\Delta = \min\{1-\beta(4+\alpha),\beta(\alpha-1)\}$.
		This concludes the proof of Theorem \ref{thm: main2}. \bigskip
	\end{proof}
	The result for the remaining values of $\alpha$, i.e. $\alpha=-(2n-1) , ~ n \in \mathbb{N}_0$, is stated in the following theorem whose detailed proof is deferred to Appendix A.
	
	\begin{thm} \label{thm:negodd}
		Assume $\alpha = -(2n-1),\, n\in \mathbb{N}_0$,  and let $a>0$. Then $\W_\alpha (z)$ is analytic in the right half-plane and there exists a polynomial $\mathbb{P} _n (z)$ of degree $2n$ and a constant $\Delta >0$ such that 
		\bb \label{negodd}
		\W _\alpha (z) = \mathbb{P}_n (z) + d_n z^{2n}\ln z + O(|z|^{2n+\Delta})
		\ee
		for $z$ small  in $V_a$, where 
		$$
		d_n = -2^{2n-1}A_n\,.
		$$
		
	\end{thm}
	
	\subsection{Coefficient asymptotics}\label{sec:3.2}
	
	Recalling that $z= \sqrt{1-4g}$, it follows from Theorems \ref{thm: main1}, \ref{thm: main2} and \ref{thm:negodd} that $\W_\alpha$ is an analytic function of $g$ in the slit-plane ${\mathbb C}\setminus [\frac{1}{4}, +\infty)$. The asymptotic behaviour of the coefficients $Z_N$ in its power series expansion \eqref{eq:Walpha1} around $g=\frac{1}{4}$, that will be needed in the next section, can be deduced by applying  transfer theorems yielding the following result.
	
	\begin{prop}
		For fixed $\alpha \in \mathbb{R}$ it holds for large $N$ that
		\bb \label{eq:asymp}
		[g^N] \W_\alpha = Z_N = C_\alpha \,N^{\frac{\alpha-3}{2}} 4^N(1+o(1))\,,
		\ee
		where the constant $C_\alpha$ is given by 
		\begin{equation} \nonumber
			C_\alpha =
			\begin{cases}
				\vspace{0.2 cm} \frac{c_\alpha}{\Gamma(\frac{\alpha-1}{2})}\,,& \text{if} ~~ \alpha>1~\text{or if}~\ -(2n+1)<\alpha<-(2n-1) ,\, n\in \mathbb{N}_0\,, \\ \vspace{0.2 cm}
				4^{n-1}|A_n|n!\,, &\text{if} ~~ \alpha=-(2n-1)\,, ~n\in \mathbb{N}_0\,. 
			\end{cases}
		\end{equation}
	\end{prop}
	
	\begin{proof} Consider first the case $\alpha \neq -(2n-1) \, , n \in \mathbb{N}_0$. With notation as in section VI.3 of \cite{flajolet2009analytic}, it follows from Theorems 
		\ref{thm: main1} and \ref{thm: main2} with $a>1$ that there exists a $\Delta$-domain 
		$$
		\Delta(\phi_a,\eta) =\{w\in\mathbb C\mid |w|<\eta, w\neq 1, |\mbox{Arg}(w-1)|>\phi_a\}
		$$
		such that 
		$$
		{\W}_\alpha(g) = \W_\alpha ^{(n)} + c_\alpha (1-4g)^{\frac{1-\alpha}{2}}\big(1+o(1)\big)\quad \mbox{for $4g\in \Delta(\eta,\phi_a)$}\,,
		$$
		where $\phi_a = \pi-2\tan^{-1} (a)<\frac{\pi}{2}$ and $\eta>0$ can be chosen arbitrarily. Applying Corollary VI.1 in \cite{flajolet2009analytic} then gives the result. 
		
		\indent For $\alpha=-2n+1$ we recall the well known fact (see e.g. section VI.2 in \cite{flajolet2009analytic}) that
		$$    
		[w^N](1-w)^n\ln\frac{1}{1-w} = (-1)^n\frac{n!}{N(N-1)(N-2)\cdots(N-n)}\,.
		$$
		This immediately implies that
		\begin{align*}
			[(4g)^N]d_n z^{2n}\ln z = \frac{2^{2n-2}(-1)^nA_n n!}{N^{n+1}} \big(1+ O(N^{-1}) \big) \,.
		\end{align*}
		Applying Theorem VI.3 of \cite{flajolet2009analytic} to the remainder term $O(|z|^{2n+\Delta})$  in \eqref{negodd}  then allows us to conclude that 
		$$
		[g^N] \W_\alpha = 2^{2n-2}(-1)^{-1}A_n n! \frac{4^N}{N^{n+1}}\big(1+o(1)\big) 
		$$
		as $N \rightarrow \infty$, as well as that $(-1)^n A_n=|A_n|$, since $Z_N$ by definition is positive. This completes the proof of the theorem.
	\end{proof}
	
	\smallskip
	
	In the particular case, $\alpha=0$, the partition function can be calculated in closed form (see e.g. \cite{durhuus2003probabilistic} for details) and is given by
	$$
	\W_0(g) = \lim_{m\to\infty}X_m(g) = \frac 12(1-\sqrt{1-4g})\,.
	$$
	Moreover, its Taylor coefficients $Z_N$ are given by the Catalan numbers,  
	\bb \label{asympCn}
	Z_N = C_{N-1} := \frac{(2(N-1))!}{N!(N-1)!} =\frac{1}{\sqrt\pi}N^{-\frac{3}{2}}4^{N-1}\big(1+O(N^{-1})\big)  \, .
	\ee
	For later use, we also note that 
	\bb\label{eq:W(0)}
	\sum_{N=1}^\infty C_{N-1} 4^{-N} = \W_0(\frac 14) = \frac 12\,.
	\ee
	
	Finally, we shall also need the asymptotic behaviour for large $N$ of the Taylor coefficients $Z_{N,M}$ of the function
	$$
	\W_{\alpha,M}(g) = \sum_{m=1}^M m^\alpha\big(X_m(g)-X_{m-1}(g)\big)\,,
	$$
	i.e.  the contribution to $\W_\alpha$ from trees of height at most $M$. Since each $X_m$ is a rational function of $g$ by \eqref{eq:sol1}, the same holds for $\W_{\alpha,M}$, and it has a unique pole closest to $g=0$ which is simple and located at  $g_M$ given by \eqref{def:gm}.  Denoting its residue by $r_M$, it follows that 
	
	\bb \label{asympZNM}
	Z_{N,M} = r_M g_M^{-(N+1)}(1+o(1))\,.
	\ee
	for $N$ large.
	
	
	\section{Local limit}

	Our goal in this section is to show the existence of the weak limit of the measures $\mu_N$ defined in \eqref{probmeasure}. 
	
	Before discussing the existence of the limit for general values of $\alpha$, let us briefly recall the familiar case $\alpha=0$. According to \eqref{probmeasure} and \eqref{asympCn} we have in this case 
	$$
	\mu_N (T) = C_{N-1} ^{-1}\,,\quad  T \in {\cal T} _N\,.
	$$
	As mentioned in the Introduction,  $\mu _N$ equals the measure obtained by conditioning a BGW tree on size $|T|=N$. More specifically, the corresponding BGW branching process is defined by the offspring distribution
	\bb\label{eq:offspring}
	p_n=2^{-n-1}\,,~n \in \mathbb{N}_0\,,
	\ee
	where $n$ denotes the number of offspring of any given individual. For details on the definition of BGW processes and basic results, the reader may consult \cite{athreyaney1972branching}, and a nice account of their local limits may be found in \cite{abraham2015introduction}. It is known \cite{aldous1998tree, kennedy1975galton} that $\mu_N $ converges weakly as $N \rightarrow \infty $ to a Borel probability measure $\mu$ concentrated on $\cal T _\infty$, and $(\cal T, \mu)$ is hence called the Uniform Infinite Planar Tree. This limiting measure $\mu$ is concentrated on single spine trees, i.e. trees consisting of a simple infinite path, the \emph{spine}, starting from the root, with finite trees (branches) grafted to its vertices on the left and right. Conditioning on the degree $k\geq 2$ of a certain spine vertex (different from the root), its spine progeny is uniformly distributed, while the $k-2$ finite branches are independent critical BGW trees with the aforementioned offspring distribution. Moreover, the degrees of spine vertices are i.i.d. random variables whose distribution equals $(k-1)\cdot 2^{-k},\,k\geq 2$.
	Alternatively, as a consequence of Theorem 2.2 in \cite{billingsley2013convergence}, the UIPT is uniquely characterised by the volume it attributes to balls in $\cT$ (see \cite{durhuus2003probabilistic} for more details). We derive an explicit expression for those volumes in Theorem \ref{thm:limit1} below. 
	
	For $\alpha\neq 0$, the measure $\mu_N$ cannot be obtained by conditioning a BGW process  on fixed size. In order to establish convergence in this more general case, we apply an approach similar to the one used in \cite{durhuus2003probabilistic} to construct the UIPT. Since $\cT$ is not compact, we proceed by showing first tightness of the sequence $(\mu_N)_{N\in\mathbb N}$, i.e. that for every $\varepsilon>0$ there exists a compact set $C\subset {\cal T}$, such that 
	\begin{equation}\label{eq:tight1}
		\mu_N({\cal T}\setminus C) < \varepsilon\quad \mbox{for all $N$}\,.
	\end{equation}
	\noindent Proving next that $(\mu_N(\cA))_{N\in\mathbb N}$ is convergent in $\mathbb R$ for any ball $\cA$ will then suffice to complete the proof. 
	
	We start by establishing lower bounds on ball volumes.
	
	\begin{prop} \label{propballvolume}
		Let $T_0 \in \mathcal T _{\rm fin}$ have height $r$ and assume it has $R$ vertices at height $r$. For every $\epsilon>0$ there exists $M_0 \in \mathbb N$, such that
		\begin{equation} \label{largeN}
			\mu_N(\mathcal{B}_{\frac 1 r } (T_0)) \geq  (1-\epsilon)R\cdot 4^{R-|T_0|} \Big( \sum\limits_{K=1} ^M C_{K-1} 4^{-K} \Big)^{R-1} \big(1+ O(N^{-1})\big)
		\end{equation}
		for all $M > M_0$.
	\end{prop}
	
\begin{figure}[ht]
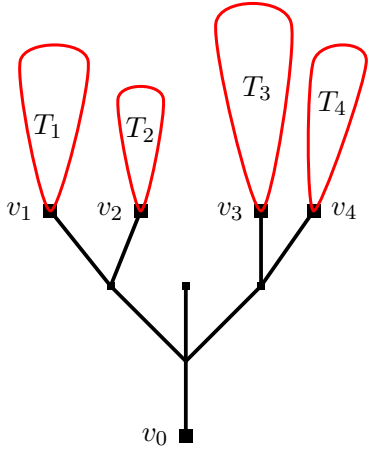

		\centering
		\begin{align*}
			\raisebox{0cm}{\treeleaf} 
		\end{align*}
		\captionof{figure}{Structure of a tree belonging to the the ball $\cB _{\frac{1}{3}}(T_0)$ around $T_0$ with root $v_0$ and leaves $(v_1,v_2,v_3,v_4)$ at height $h(T_0)=3$, obtained by grafting 4 trees $T_1,T_2,T_3,T_4 \in \cT$  onto $T_0$ at $v_1,v_2,v_3,v_4$, respectively.}
		\label{figure:treeleaf}
	\end{figure}

	\begin{proof}
		For simplicity, let us denote $\mathcal{B}_{\frac{1}{r}}(T_0)$ by $\cA$, and let $v_1,...,v_R$ denote the vertices of $T_0$ at height $r$. Then any $T \in \cA$ can be obtained by grafting a sequence $(T_1,...,T_R) \in \mathcal{T}^R$ of trees onto $T_0$ such that the root edge of $T_i$, $1\leq i \leq R$, is identified with the edge in $T_0$ incident on $v_i$, see Figure \ref{figure:treeleaf} for an illustration. Call $T_1, \dots, T_R$ the branches of $T$. Given $M>0$ and $i_0\in\{1,\dots,R\}$, let the set ${\cA}^M_{i_0}$ be given by
		$$
		{\cA}^M_{i_0} =\{T\in \cA\mid |T_i|\leq M\, \mbox{for}\, i\neq i_0\,\;\mbox{and}\; h(T_{i_0}) > M\}\,.
		$$
		For $N>|T_0|+RM$, the set  $A^M_{i_0}\cap {\cal T}_N$ is non-empty and consists of trees whose highest and largest (in terms of size) branch is $T_{i_0}$. Moreover, given $\epsilon>0$, we can choose $M_0$ sufficiently large such that for $M>M_0$
		\bb\label{eq:ineqh}
		\Big(\frac{h(T)}{h(T_{i_0})}\Big)^\alpha  = \Big(\frac{h(T_{i_0}) + r-1}{h(T_{i_0})}\Big)^\alpha\geq 1-\epsilon \quad\mbox{for all $T\in \cA^M_{i_0}$}.
		\ee
		\noindent Clearly, the sets ${\cA}_i^M,\, i=1,\dots,R$, are disjoint subsets of $\cA$, if $N> |T_0|+RM$, and hence
		\bb\label{eq:ineqmu2}
		\mu_N(\cA) \geq \sum_{i_0 = 1}^R \mu_N({\cA}^M_{i_0})\,.
		\ee 
		Moreover, for such $N$ 
		\bb \label{muA0}
		\mu_N({\cA}^M_{i_0}) =  \sum_{\substack{N_1+\dots+N_R=N+R-|T_0|\\ N_i\leq M,\,i\neq i_0}} \frac{\tilde Z_{N_{i_0},M}}{Z_N} \prod_{i \neq i_0} C_{N_i-1}\,,
		\ee
		where the modified partition function $\tilde Z_{K,M}$ is given by 
		\bb \label{eq:tildeZ}
		\tilde Z_{K,M} = \sum _{\substack{T \in {\mathcal T}_K\\h(T)> M}} (h(T)+r-1)^\alpha \geq (1-\epsilon)\sum _{\substack{T \in {\mathcal T}_K\\h(T)> M}} h(T)^\alpha\,.
		\ee
		Here, the last inequality is a consequence of \eqref{eq:ineqh}, assuming that $M>M_0$. With notation as in section \ref{sec:3.2}, the last sum in \eqref{eq:tildeZ} equals $Z_K-Z_{K,M}$, and according to  \eqref{eq:asymp} and \eqref{asympZNM} the bound
		\bb\label{intrat}
		\frac{Z_{N_{i_0},M}}{Z_N} \leq  \mbox{\rm cst}\cdot \frac{N^{\tfrac{3-\alpha}{2}}}{4^Ng_{M}^N}
		\ee
		holds for fixed $M$ and for $N$ large enough, where we have also taken into account that $g_M\leq 1$ and $N_{i_0} < N$.
		Since $g_M>\frac 14$, the right-hand side of \eqref{intrat} decays exponentially with $N$. It thus follows from \eqref{eq:asymp}, for fixed $M$ and fixed $N_{i},i\neq i_0$, fulfilling the summation constraints of \eqref{muA0}, that 
		\bb \nonumber 
		\frac{\tilde Z_{N_{i_0},M}}{Z_N}  =  4^{-\sum\limits_{i\neq i_0} N_i + R-|T_0|}\big(1+ O(N^{-1})\big)\,.
		\ee   
		Using this together with \eqref{eq:tildeZ} and \eqref{eq:ineqmu2}, we get for any $M>M_0$ and $N$ sufficiently large that 
		\bb \nonumber \label{muA1}
		\mu_{N}(\cA) \geq (1-\epsilon) \cdot 4^{R-|T_0|}  \sum_{i_0=1}^R\sum_{\substack{1\leq N_i\leq M\\ i\neq i_0}}  \prod_{i \neq i_0} 4^{-N_i} C_{N_i-1}\big(1+O(N^{-1})\big)\,,
		\ee 
		which is equivalent to \eqref{largeN}.
	\end{proof}
	
	For $\epsilon>0$ and $M>M_0$, let us denote the large-$N$ limit of the right hand side of \eqref{largeN} by $\Lambda(T_0,M,\epsilon)$, i.e.
	\bb \label{lambdaeps}
	\Lambda(T_0,M, \epsilon) = (1-\epsilon)R\cdot 4^{R-|T_0|} \Big( \sum\limits_{K=1} ^M C_K 4^{-K} \Big)^{R-1} \,,
	\ee
	and let
	\bb \label{defLambda}
	\Lambda(T_0) := \lim\limits_{M \rightarrow \infty} \Lambda(T_0,M,0) = R\cdot 2^{R+1} 4^{-|T_0|} \, .
	\ee
	
	\begin{lem}\label{lem:addupto} 
		For any $r\in\mathbb N$ it holds that 
		\bb\label{addupto}
		\sum_{\substack{T_0 \in \cT _{\text{\rm fin}}: h(T_0)=r}} \Lambda(T_0)  \; =\; 1\,.
		\ee
	\end{lem}
	
	\begin{proof}
		We use an inductive argument. For $r=1$ the statement trivially holds, 
		so let $r \geq 2$ be arbitrary and assume \eqref{addupto} holds for $r-1$. Recall that the $R$-factor in \eqref{defLambda} originates from summing over the position $i_0$ of the long branch out of $R$ branches. This branch has an ancestor $j$ in $T_0$ at height $r-1$, which is the root of the branch $T_{i_0}$. Letting $R'$ denote the number of vertices at height $r-1$, it follows that summing over trees $T_0$ of height $r$ with $T_0 ^\prime := B_{r-1}(T_0)$ fixed and with a marked root edge of the large branch amounts to summing over all possible choices of the remaining $R-1$ edges at maximal height. Since the number of such choices equals $\binom{R+R^\prime-1}{R^\prime}$, we get
		\begin{align} \nonumber
			\sum_{\substack{T_0 : B_{r-1}(T_0) = T_0 ^\prime\\ h(T_0)=r}} \Lambda(T_0) 
			&=  4 ^{-|T_0^\prime|}R^\prime \sum_{R \geq 1} \binom{R+R^\prime-1}{R^\prime} (\half)^{R-1} \\ 
			&= R^\prime \cdot 2^{R^\prime +1} 4^{-|T_0 ^\prime |}  \,.
			\label{identical}
		\end{align}
		Since the last expression in \eqref{identical} equals $\Lambda(T_0^\prime)$, this completes the proof.
	\end{proof}	
	
	\begin{cor}\label{tightness0} The sequence  $(\mu_N)_{N\in\mathbb N}$ of measures on $\cT$ given by \eqref{probmeasure} is tight.
	\end{cor}
	
	\begin{proof}  
		It is easy to verify (see e.g. \cite{durhuus2003probabilistic} for details) that sets of the form 
		\bb \label{Compact}
		C = \bigcap _{r=1} ^\infty \{ T \in \mathcal T \Big| ~ |B_r (T)|  \leq K_r \}\,,
		\ee
		where $(K_r)_{r\in\mathbb N}$ is any sequence of positive numbers, are compact. In order to establish \eqref{eq:tight1} for sets of this form, it is sufficient to show that for every $\delta > 0$  and $r \in \mathbb{N}$ there exists $K_r > 0 $, such that
		\bb\label{Ddelta}
		\mu_N (\{ T \in \mathcal T \Big| ~ |B_r(T)| > K_r \} ) < \delta\,,\quad N \in \mathbb{N}\,.
		\ee
		Indeed, choosing $\delta$ to be $r$-dependent of the form $\delta_r=\frac{\varepsilon}{2^r}$ and defining $C$ by \eqref{Compact} for the corresponding values of $K_r$ determined by \eqref{Ddelta}, we obtain  
		\bb\label{Cepsilonbound}
		\mu_N (\cT\setminus C)\;\leq\; \sum_{r=1}^\infty \mu_N (\{T\in\cT \mid |B_r(T)|> K_r\}) \leq \epsilon\,.
		\ee
		
		To establish \eqref{Ddelta}, let $r\geq 1$ be given and choose first by Lemma \ref{lem:addupto} a finite subset $\cT _0$ of $\cT _{\rm fin}$ consisting of trees of height $r$,  such that
		\bb \nonumber
		\sum_{T_0\in \cT_0} \Lambda(T_0) \geq 1-\delta\,.
		\ee
		Then apply \eqref{lambdaeps} and Proposition \ref{propballvolume} to find $\epsilon>0$ small enough and $M$ large enough, such that
		\bb \nonumber
		\sum_{T_0\in \cT_0} \Lambda(T_0,M,\epsilon) \geq 1-2\delta\,.
		\ee
		Using \eqref{largeN}, we can then choose $N_0$  such that 
		\bb\label{1minus3delta}
		\sum_{T_0\in\cT_0}  \mu_N (\mathcal{B} _{\frac{1}{r}}(T_0)) \geq 1-3\delta \quad \mbox{for $N\geq N_0$}\,.
		\ee
		Now, let $K_0=\mbox{max}\{|T_0|\mid T_0\in\cT_0\}$ and observe that the pairwise disjoint balls $\mathcal{B} _{\frac{1}{r}}(T_0),\, T_0\in\cT_0$, are contained in $\{T\in\cT \mid |B_r(T)|\leq  K_0\}$. Therefore,  \eqref{1minus3delta} implies
		$$
		\sum_{K=1}^{K_0} \mu_N(\{T\in\cT \mid |B_r(T)| =K\}) \geq 1-3\delta \quad \mbox{for $N\geq N_0$}\,.
		$$
		Since the sets $\{T\in\cT \mid |B_r(T)| =K\},\, K\in\mathbb N$, are pairwise disjoint, it follows that 
		$$
		\mu_N(\{T\in\cT \mid |B_r(T)| > K\}) \leq 3\delta\quad  \mbox{for $N\geq N_0$ and $K>K_0$}.
		$$
		Choosing $K>N_0$, this inequality holds for all $N$, and thus the proof is complete.
	\end{proof}
	\begin{thm}\label{thm:limit1}
		The sequence  $(\mu _N)$ defined by \eqref{probmeasure} is weakly convergent to a Borel probability measure $\mu$ on $\cT$, that is characterised by 
		\bb\label{limit1}
		\mu(\mathcal{B}_{\frac 1r}(T_0)) = \Lambda(T_0) = R\cdot 2^{R+1} 4^{-|T_0|}\,,
		\ee
		for any  tree $T_0\in\cT_{\rm fin}$ of height $r$, where $R$ denotes the number of vertices at height $r$. Moreover, $\mu$ is equal to the UIPT.
	\end{thm}

	\begin{proof}
		Since $(\mu _N)$ is tight, by the previous corollary, it has a weakly convergent subsequence $(\mu _{N_i})$ converging to a probability measure $\mu $ on $\cT$. We shall show that the limit $\mu$ is independent of the subsequence and hence that $(\mu _N)$ is convergent. Since the balls in $\cT$ have empty boundary, it follows from Theorem 2.1  in \cite{billingsley2013convergence} that $\mu _{N_i}(\mathcal{B}_{\frac 1r}(T_0))$ converges to $\mu(\mathcal{B}_{\frac 1r}(T_0))$ as $i\to\infty$. Using Proposition \ref{propballvolume}, this implies 
		\bb \nonumber
		\mu(\mathcal{B}_{\frac{1}{r}}(T_0) )\geq \Lambda(T_0,M,\epsilon)
		\ee
		for any $T_0\in\cT_{\rm fin}$ of height $r \in \mathbb N$, and any $\epsilon>0$ and $M >M_0$. Letting $M\to\infty$ and subsequently $\epsilon\to 0$, we obtain 
		\bb \label{ineqeq}
		\mu(\mathcal{B}_{\frac{1}{r}}(T_0) )\geq \Lambda(T_0)\,.
		\ee
		Finally, using Lemma \ref{lem:addupto} and the fact that $\mu$ is a probability measure, it follows that equality holds in \eqref{ineqeq}. Since any Borel probability measure on $\cT$ is uniquely determined by its value on balls, by Theorem  2.2 in \cite{billingsley2013convergence}, we have shown that the limit $\mu$ is unique and independent of $\alpha$. In particular, $\mu$  equals the UIPT, and the proof is complete.
	\end{proof}
	
	\section{Concluding remarks}
	We note that, despite the widely varying singular behaviour of the generating function $\W_\alpha (g)$: it is finite at the critical point for $\alpha<1$, has a logarithmic divergence when $\alpha=1$ and a power-like divergence for $\alpha>1$, we have found that the local limit of the distributions \eqref{probmeasure} is independent of the exponent $\alpha \in \mathbb{R}$. Whether more general BGW trees respond in a similar way to a powerlike height coupling, we do not know, since our approach relies on knowing the explicit form of the generating function for fixed-height partition functions, given in \eqref{eq:sol2}, the analogue of which is not generally available. 
	It is, however, conceivable that more general techniques based on recursion relations alone could be developed.
	
	Another way of extending the results of this paper would be allowing different forms of height couplings. In \cite{meltem2021height}, we consider weights of exponential form, $k^h$, at fixed size, partly motivated by findings in \cite{durhuus2021critical}, where an analysis of certain statistical mechanical models of loops on random so-called causal triangulations is carried out. Via a bijective correspondence between causal triangulations and rooted planar trees, it turns out that some of those models can be related to models of planar random trees with exponential height coupling with $k>1$. As shown in \cite{meltem2021height}, it turns out that the local limits exhibit qualitatively different behaviours, depending on whether $0<k<1,\, k=1$, or $k>1$. 
	
	\bigskip
	
	\subsection*{Acknowledgements}
	
	The authors acknowledge support from Villum Fonden via the QMATH Centre of Excellence (Grant no.~10059).
	
	\bigskip
	
	\section*{Appendix}
	\emph{Proof of Theorem \ref{thm:negodd}.}
	We shall make use of the fact that  
	\bb \label{claimint}
	\int \limits_r ^\infty \Big(\frac{1}{t^{2n-1}(\cosh(2t)-1)} - \frac{1}{t^{2n-1}}L_{n-1}(2t)\Big) dt = \zeta -2^{2n-1}A_n \ln r + O(r^2)
	\ee
	for small $r>0$, where $\zeta$ is a constant and $A_n$ and $L_n$ are defined in \eqref{laurent} and \eqref{eq:Ln}, respectively. Clearly, the integral in \eqref{claimint} is convergent for $r>0$, and since the Laurent expansion \eqref{laurent} is convergent for $0<|t|<\pi$ we get that, up to an additive constant, the integral equals
	$$
	\int \limits_r ^{\pi/2} \sum_{k=n}^\infty A_k 2^{2k-1} t^{2(k-n)-1}dt = - 2^{2n-1}A_n\ln r+ C + O(r^2)\,,
	$$
	where $C$ is a constant, and thus \eqref{claimint} follows. On the other hand, using Cauchy's theorem as previously, we can rewrite the integral in \eqref{claimint} as a line integral along the circular arc $C_r$  of radius $r$ centered at $0$ connecting $r$ and $r\frac{z}{|z|}$ and along the half line $\ell_z: s \rightarrow sz$ with endpoint at $r\frac{z}{|z|}$ inside the wedge $V_a$, and get 
	\begin{align} \nonumber
		\int \limits_{\frac{r}{|z|}}^\infty \Big(\frac{z}{(sz)^{2n-1}(\cosh(2sz)-1)}&-\frac{z}{(sz)^{2n-1}}L_{n-1}(2sz) \Big)ds\\ \label{cauchy2}
		& = \zeta -2^{2n-1}A_n (\ln r +i\,\text{Arg}\,z) + O(r^2)\,,
	\end{align}
	where we have used that
	$$
	\int  \limits_{C_r} \Big(\frac{1}{w^{2n-1}(\cosh(2w)-1)} -\frac{1}{w^{2n-1}} L_{n-1}(2w)\Big)dw  = i\,2^{2n-1} A_n\,\text{Arg} z + O(r^2) \,.
	$$ 
	
	Recalling the definition \eqref{Wn} of $\W_\alpha^{(n)}$, we have 
	$$
	\W_\alpha(z) - \W_\alpha^{(n-1)}(z) = \sum_{m=1}^\infty m^\alpha\Big(\frac{z^2}{D_m(z)}-T_{n-1}^m(z)\Big)\,,
	$$
	where $T_{n-1}^m$ is defined by \eqref{eq:Tnm}. Consider first the contribution $S_1$ to this sum from $m\leq \frac{r}{|z|}$, where $r$ satisfies \eqref{eq:cond1}. With notation as in  Lemma \ref{lem: lemma2} b) we rewrite
	\begin{align*}
		S_1 &= \sum \limits_{m|z|\leq r} \frac{1}{2m^{2n}(m+1)} \sum \limits_{k=n} ^\infty (c^m _{2k})(2mz)^{2k} \\
		&= \sum \limits_{m|z|\leq r}  \Big[ 2^{2n-1}z^{2n}\Big(\frac{A_n}{m+1} + \frac{c_n ^m-A_n}{m+1}\Big) + \sum \limits_{k=n+1} ^\infty \frac{c^m _{2k}}{2m^{2n}(m+1)}(2mz)^{2k} \Big]\,.
	\end{align*}
	Here, the first term inside round parenthesis is harmonic and yields the contribution 
	$$
	2^{2n-1} z^{2n} A_n \Big(\ln\frac{r}{|z|}-1+\gamma) + O\big(\frac{|z|}{r}\big)\Big)\,,
	$$
	while the bound \eqref{eq:c_2k2} implies that the contribution from the second term equals
	$$
	e_n z^{2n} + O\big(\frac{|z|^{2n+1}}{r}\big)\,, 
	$$
	where 
	$$
	e_n =   2^{2n-1}\sum_{m=1}^\infty\frac{c_n ^m-A_n}{m+1}\,.
	$$
	Using the bound \eqref{eq:lem3b}, the remaining contribution to $S_1$ is easily bounded by  $O(r^2|z|^{2n}) $, such that we obtain
	\bb \nonumber  \label{S11}
	S_1 =  2^{2n-1} z^{2n} A_n (\ln\frac{r}{|z|} -1+\gamma) +  e_n z^{2n} + O(r^2|z|^{2n}) +  O\big(\frac{|z|^{2n+1}}{r}\big)\,.
	\ee
	Taking into account \eqref{cauchy2}, we conclude that 
	
	\begin{align} \nonumber  \label{2piece2}
		\W_\alpha (z) -\W_\alpha^{(n-1)}(z) - (e_n + \zeta + 2^{2n-1}&A_n(\gamma-1)) z^{2n} + 2^{2n-1} A_n z^{2n}\ln z \\ \nonumber 
		&= 		S - I  +  O(r^2|z|^{2n}) +  O\big(\frac{|z|^{2n+1}}{r}\big)\,,
	\end{align}
	where
	$$
	S =  \sum_{m|z|> r} \Big(\frac{z^2}{m^{2n-1}D_m(z)}-\frac{1}{m^{2n-1}}T_{n-1}^m(z)\Big)
	$$
	and
	$$
	I = \int \limits_{\frac{r}{|z|}} ^\infty \Big(\frac{z^2}{s^{2n-1}(\cosh(2sz)-1)}-\frac{z^2}{s^{2n-1}}L_{n-1}(2sz)\Big)ds\,.
	$$
	
	Next, we proceed to estimate $|S-I|$ as before by splitting the summation and integration domains corresponding to $r < m|z|\leq \frac 1r$ and $m|z| > \frac 1r$ and similarly for $s$. Calling the corresponding sums and integrals $S_2, S_3$ and $I_2, I_3$, respectively, we first rewrite   
	$$
	S_2-I_2 = A + B + C\,,
	$$
	where 
	\begin{align*}
		A&:= \sum \limits_{r<m|z|\leq r^{-1}} \frac{z^2}{m^{2n-1}} \Big(\frac{1}{D_m(z)} - \frac{1}{\cosh(2mz)-1} \Big) \\
		B&:= \sum \limits_{r<m|z|\leq r^{-1}} \frac{z^2}{m^{2n-1}(\cosh(2mz)-1)} - \int \limits _{\frac{r}{|z|}} ^{\frac{1}{r|z|}} \frac{z^{2}ds}{s^{2n-1}(\cosh(2sz)-1)} \\
		C&:=  \int \limits_\frac{r}{|z|}^{\frac{1}{r|z|}} \frac{1}{2 s^{2n+1}} \Big[ 1+\sum_{k=1} ^{n-1} A_k (2sz)^{2k} \Big]ds -\sum \limits_{r<m|z|\leq r^{-1}} \frac{1}{2m^{2n}(m+1)}\Big[1+ \sum \limits _{k=1} ^{n-1} c^m _{2k} (2mz)^{2k}\Big]\,.
	\end{align*}
	Thus, $A$ and $B$ fulfill the bounds \eqref{eq:estA} and \eqref{eq:estB}, valid for $\alpha\leq1$, i.e.
	\bb \label{An}
	|A| \leq \text{cst} \cdot \frac{|z|^{2n+1}}{r^{2n+3}} \,, \quad |B| \leq \text{cst} \cdot \frac{|z|^{2n+1}}{r^{2n+3}} \,,
	\ee
	provided that \eqref{eq:cond2} is fulfilled.
	Similarly the bound \eqref{eq:estC} holds for $C$, 
	\bb \label{Cn}
	|C| \leq \text{cst} \cdot \frac{|z|^{2n+1}}{r^{2n+1}}\,,
	\ee
	provided that \eqref{cond3} is fulfilled. 
	
	For $S_3 - I_3$, the bounds \eqref{S31},\eqref{I31},\eqref{Wn2} and \eqref{Tn} are still applicable. Noting that the index $n$ of $T_n^m$ and         $L_n$  in the latter two is replaced by $n-1$, one obtains 
	\bb \label{s3s3ni3}
	|S_3  -I_3| \leq \text{cst} \cdot |z|^{2n}\big(\exp(-\frac{\text{cst}}{r}) + r^2\big) \,.
	\ee
	
	Collecting \eqref{An}, \eqref{Cn} and \eqref{s3s3ni3} and defining 
	$$
	{\mathbb P}_n(z) = \W_\alpha^{(n-1)}(z) + 2^{2n-1}(e_n+\zeta + (A_n(\gamma-1))z^{2n}\,,
	$$
	we conclude that 
	\bb 
	\big| \W_\alpha (z)- {\mathbb P}_n(z) +2^{2n-1}A_n z^{2n} \ln z\big|
	\,\leq\, \text{cst} \cdot \Big(r^2 |z|^{2n} + \frac{|z|^{2n+1}}{r^{2n+3}}\Big)\,, \label{snfinal}
	\ee
	provided $z$, $r$ and $n$ fulfill \eqref{eq:cond1}, \eqref{eq:cond2} and \eqref{cond3}.  
	Choosing $r=|z|^\beta$, where $0<\beta<1$, these conditions are satisfied for $|z|$ small enough. Finally, imposing in addition 
	$$
	\beta(2n+3)<1\,,
	$$
	it follows from \eqref{snfinal} that \eqref{negodd} holds with $\Delta = \min\{1-\beta(2n+3),2\beta \}$.
	
	This concludes the proof of Theorem \ref{thm:negodd}. \qed

	\bibliographystyle{abbrv}

\end{document}